
\input amstex.tex
\documentstyle{amsppt}
\magnification=\magstep1
\hsize=12.5cm
\vsize=18cm
\hoffset=1cm
\voffset=2cm

\font\ff=cmr9

\def\DJ{\leavevmode\setbox0=\hbox{D}\kern0pt\rlap
{\kern.04em\raise.188\ht0\hbox{-}}D}
\footline={\hss{\vbox to 2cm{\vfil\hbox{\rm\folio}}}\hss}
\nopagenumbers
\font\ff=cmr8
\def\txt#1{{\textstyle{#1}}}
\baselineskip=13pt
\def\hf{{\textstyle{1\over2}}}
\def\a{\alpha}\def\b{\beta}
\def\d{{\,\roman d}}
\def\e{\varepsilon}

\def\g{\gamma} 

\def\s{\sigma}

\def\={\;=\;}

\def\D{\Delta}
\def\no{\noindent}
  \def\s{\sigma}

\def\no{\noindent}
\def\e{\varepsilon}
\def\D{\Delta}
\def\no{\noindent}
\def\e{\varepsilon}

\def\no{\noindent}
\font\teneufm=eufm10
\font\seveneufm=eufm7
\font\fiveeufm=eufm5
\newfam\eufmfam
\textfont\eufmfam=\teneufm
\scriptfont\eufmfam=\seveneufm
\scriptscriptfont\eufmfam=\fiveeufm
\def\mathfrak#1{{\fam\eufmfam\relax#1}}
\define\pahx{\frac{4\pi}{k}\sqrt{nx}}
\define\pahax{(\pahx)}

\font\tenmsb=msbm10
\font\sevenmsb=msbm7
\font\fivemsb=msbm5
\newfam\msbfam
\textfont\msbfam=\tenmsb
\scriptfont\msbfam=\sevenmsb
\scriptscriptfont\msbfam=\fivemsb
\def\Bbb#1{{\fam\msbfam #1}}

\def \NN {\Bbb N}
\def \CC {\Bbb C}

\def\rightheadline{{\hfil{\ff
The Laplace transform of the square in the circle problem}\hfil\tenrm\folio}}

\def\leftheadline{{\tenrm\folio\hfil{\ff
Aleksandar Ivi\'c }\hfil}}
\def\emptyheadline{\hfil}
\headline{\ifnum\pageno=1 \emptyheadline\else
\ifodd\pageno \rightheadline \else \leftheadline\fi\fi}

\topmatter
\title A NOTE ON THE LAPLACE TRANSFORM OF
THE SQUARE IN THE CIRCLE PROBLEM \endtitle
\author   Aleksandar Ivi\'c \endauthor
\address{
Aleksandar Ivi\'c, Katedra Matematike RGF-a
Universiteta u Beogradu, \DJ u\v sina 7, 11000 Beograd,
Serbia (Yugoslavia).}
\endaddress
\keywords Circle problem, Laplace transform,
additive problems
\endkeywords
\subjclass Primary 11N37; Secondary 44A10 \endsubjclass
\email {\tt aleks\@ivic.matf.bg.ac.yu,
aivic\@rgf.rgf.bg.ac.yu} \endemail
\abstract
{If $P(x)$ is the error term in the circle problem, then it is proved
that
$$
\int_0^\infty P^2(x)e^{-x/T}\d x ={1\over4}\left({T\over\pi}\right)
^{3/2}\sum_{n=1}^\infty r^2(n)n^{-3/2} - T + O_\e(T^{\frac{2}{3}+\e}),
$$
improving the earlier result with exponent $\frac{5}{6}$ in the error term.
The new bound is obtained by using  results of F. Chamizo on the
correlated sum $\sum_{n\le x}r(n)r(n+h)$, where $r(n)$
is the number of representations of $n$ as a sum of
two integer squares.}
\endabstract
\endtopmatter

\vglue 1cm
\head 1. Introduction
\endhead

\bigskip\no
Let $r(n) = \sum_{n=a^2+b^2}1$ denote the number of representations of
$n\;(\in\NN)$ as a sum of two integer squares. Thus ${1\over4}r(n)$ is
multiplicative and
$$
r(n) \;=\;  4\sum_{d|n}\chi(d),   \leqno(1.1)
$$
where $\chi(n)$ is the non-principal character mod$\,4$. A classical
problem, with a rich history, is the circle
problem. It consists of the estimation of the function
$$
P(x) \;=\; {\sum_{n\le x}}'r(n) - \pi x + 1,\leqno(1.2)
$$
where, as usual, ${\sum\limits_{n\le x}}'$ means that the last term
in the sum is to be halved if $x$ is an integer. One can estimate
$P(x)$ pointwise and in the mean square sense. M.N. Huxley [3] proved
that
$$
P(x) \= O(x^{23/73}\log^cx)\qquad(c > 0,\; \frac{23}{73} = 0,3150684\ldots\,),
\leqno(1.3)
$$
which is the last in a series of improvements by the estimation
of intricate exponential sums. The mean square formula for $P(x)$ is
written in the form
$$
\int_0^XP^2(x)\d x = \left(\frac{1}{3\pi^2}\sum_{n=1}^\infty r^2(n)n^{-3/2}
\right)X^{3/2} + Q(X),\leqno(1.4)
$$
where $Q(X)$ is considered as the error term. The best known bound is
$$
Q(X) \;=\;O(X\log^2X),\leqno(1.5)
$$
proved long ago by I. K\'atai [9]. From (1.4) and (1.5) one deduces
that
$$
P(X) \;=\;\Omega(X^{1/4}),\leqno(1.6)
$$
where as usual $f = \Omega(g)$ means that $\lim_{x\to\infty}f(x)/g(x)
 \not= 0$. The omega-result (1.6) favours the long standing conjecture
that
$$
P(X) \;=\;O_\e(X^{{1\over4}+\e}),\leqno(1.7)
$$
where $\e$ denotes arbitratily small positive numbers, not necessarily
the same ones at each occurrence. A comparison of (1.3) and (1.7) shows
that there is a big gap between the known and conjectured pointwise
estimates for $P(x)$.

\medskip
A useful representation of $P(x)$ is the classical formula
$$
P(x) \= x^{1/2}\sum_{n=1}^\infty r(n)n^{-1/2}J_1(2\pi\sqrt{xn}),
\leqno(1.8)
$$
due to G.H. Hardy [2], where $J_1$ is the customary Bessel function.
The series in (1.8) is boundedly, but not absolutely convergent.
This causes problems in practice, and one can use the truncated form
$$
P(x) \= -\frac{x^{1/4}}{\pi}\sum_{n\le N}r(n)n^{-3/4}
\cos(2\pi\sqrt{xn} + {\pi\over4}) + O_\e(x^\e + x^{{1\over2}+\e}N^{-
{1\over2}}),\leqno(1.9)
$$
which is valid for $x \ge 2,\, 2 \le N \le x^A,$
and $A > 0$ is any constant. Trivial estimation of the sum in (1.9)
(with $N = x^{1/3}$) yields at once the bound $P(x) \ll_\e x^{{1\over3}+\e}$.
\no
\bigskip
{\bf Acknowledgement.} I wish to thank F. Chamizo and T. Meurman for
valuable remarks.
\bigskip

\head 2. The Laplace transform of $P^2(x)$
\endhead

\bigskip\no
The difficulties encountered in evaluating mean square integrals like
the one in (1.4) are less pronounced    when the integrand is
multiplied by an appropriate smooth function. In [5] the Laplace
transforms of $P^2(x)$ and $\D^2(x)$ were evaluated, when $s = 1/T\to 0+$,
and
$$
\D(x) = {\sum_{n\le x}}'d(n) - x(\log x +2\g -1) - {\txt{1\over4}},
\quad d(n) = \sum_{\delta|n}1\leqno(2.1)
$$
is the error term in the classical Dirichlet divisor problem ($\g$
is Euler's constant). It was proved that
$$
\int_0^\infty P^2(x)e^{-x/T}\d x ={1\over4}\left({T\over\pi}\right)
^{3/2}\sum_{n=1}^\infty r^2(n)n^{-3/2} - T + O_\e(T^{\a+\e})\leqno(2.2)
$$
and
$$\eqalign{
&\int_0^\infty \D^2(x)e^{-x/T}\d x =\cr&
= {1\over8}\left({T\over2\pi}\right)
^{3/2}\sum_{n=1}^\infty d^2(n)n^{-3/2}
+ T(A_1\log^2T + A_2\log T + A_3) + O_\e(T^{\b+\e}).\cr}\leqno(2.3)
$$
The $A_j$'s are suitable constants ($A_1 =
-1/(4\pi^2)$), and the constants $\hf \le \a < 1$
and $\hf \le \b < 1$ are defined by the asymptotic formula
$$
\sum_{n\le x}r(n)r(n+h) = \frac{(-1)^h8x}{h}\sum_{d|h}(-1)^dd
+ E(x,h), E(x,h) \ll_\e x^{\a+\e},\leqno(2.4)
$$
$$
\sum_{n\le x}d(n)d(n+h) = x\sum_{i=0}^2(\log x)^i\sum_{j=0}^2c_{ij}
\sum_{d|h}\left({\log d\over d}\right)^j + D(x,h),
D(x,h) \ll_\e x^{\b+\e}.\leqno(2.5)
$$
The $c_{ij}$'s are certain absolute constants, and the  $\ll$--bounds
both in (2.4) and in (2.5) should hold uniformly in $h$ for $1\le h
\le x^{1/2}$. With the values $\a = 5/6$ of D. Ismoilov [6] and
$\b = 2/3$ of Y. Motohashi [11] it followed then that (2.2) and (2.3)
hold with $\a = 5/6$ and $\b = 2/3$. Motohashi's fundamental paper
(op. cit.)  used the powerful methods of spectral theory of the
non-Euclidean Laplacian. A variant of this approach was used recently
by T. Meurman [10] to sharpen Motohashi's bound for $D(x,h)$ for
`large' $h$, specifically for $x^{7/6} \le h \le x^{2-\e}$, but the
limit of both methods is $\b = 2/3$ in (2.5).

Although one expects, by general analogies between the circle and
divisor problems (see e.g., [4, Chapter 13]), that $\a = \b$ holds
(and that in fact $\a =  \b = \hf$), proving this is difficult.
If one wants to generalize the method of Motohashi or Meurman to
$E(x,h)$, one encounters several difficulties. One stems from the fact
that $r(n)$ is given by (1.1), while $d(n) = \sum_{\delta|n}1$
contains no characters. This is reflected in the following. Namely
Meurman uses a Voronoi--type formula for sums of $d(n)F(n)$
($F(x) \in C^1[a,b]$) when $n$ lies in a given residue class. Such
a formula is easily derived from the summation formula (see M. Jutila [7])
$$
{\sum_{a\le n\le b}}'d(n)e(\frac{nh}k)F(n)
=\frac1k\int_a^b (\log x +2\gamma -2\log k)F(x)\d x
$$
$$+\frac1k\sum_{n=1}^\infty d(n)\int_a^b
\bigl(-2\pi e(\frac{-n\bar{h}}{k}) Y_0\pahax
+4e(\frac{n\bar{h}}{k})K_0\pahax \bigr)F(x)\d x,
$$
which is valid for $0 < a < b,\,F(x) \in C^1[a,b]$ and $(h,k) = 1$.
However, the analogue of this formula for sums of $r(n)e(\frac{nh}k)F(n)$
is not so simple arithmetically. Namely  M. Jutila analyzed
this problem in his paper [8]. His equations (27) and (28) give
$$\eqalign{
{\sum_{a\le n\le b}}'r(n)e({nh\over k})F(n) &= \pi k^{-2}G_Q(k,h)
\int_a^b F(x)\d x\cr&
+ 2\pi (2k)^{-1}\sum_{k=1}^\infty \tilde{r}(n)e(-\overline{4h}{n\over k})
\int_a^b J_0({2\pi\over k}\sqrt{xn})F(x)\d x,\cr}
$$
where $\overline{h}$ is the multiplicative inverse of $h$ mod$\,k$ and
$$
G_Q(k,h) \= \left(\sum_{x=1}^ke({h\over k}x^2)\right)^2
$$
is the square of the Gauss sum, so it is zero for $k = 4m+2$ and
$\chi(k)k$ for $k = 4m +1$. When $k = 1$ we do get the `ordinary'
Voronoi formula for $r(n)$ (in which case $k^{-2}G_Q(k,h) = 1$),
but for general $k$ the function $\tilde{r}(n)$ (it is small, being
$\le 2r(n) \ll_\e n^\e$) depends also on $k$. The outcome of this
summation formula will be that  we shall  not
get the `nice' Kloosterman sum as happened in the case of $d(n)$, but
some `twisted'  sums. In the case of $d(n)$ one used Kuznetsov's trace
formula for sums of Kloosterman sums, but in the case of $r(n)$ the
analogue of this step is hard.

\medskip
Nevertheless we can avoid these difficulties and appeal to results
of F. Chamizo [1]  to show that $\a = 2/3$ is indeed possible in
(2.2), which is the limit of present methods coming from the use
of spectral theory. Thus we have the following

\medskip
{\bf THEOREM}. {\it We have}
$$
\int_0^\infty P^2(x)e^{-x/T}\d x ={1\over4}\left({T\over\pi}\right)
^{3/2}\sum_{n=1}^\infty r^2(n)n^{-3/2} - T + O_\e(T^{{2\over3}+\e}).
\leqno(2.6)
$$

\bigskip\no
\head 3. Proof of the Theorem
\endhead
We shall follow the method of [5] and use  Theorem
4.3 of F. Chamizo [1]. This says that, uniformly
for arbitrary $\a_m \in \CC$ and $M > 1, \,N > 1$,
$$
\sum_{M<m\le2M}\a_mE(N,m) \ll_\e ||\a||_2(N^{2/3+\e}M^{1/2}
+ N^{1/3}M^{5/6+\e}),\leqno(3.1)
$$
where $||\a||_2 = \left(\sum_{M<m\le2M}|\a_m|^2\right)^{1/2}$ is the norm
of the sequence $\{\a_m\}$. We also have by [1, Theorem 5.2] the pointwise
estimate
$$
E(N,m) \ll_\e N^{{2\over3}+\e}m^{5\over42}\qquad(m \le N).\leqno(3.2)
$$
Actually Chamizo defines (see (2.4))
$$
E(N,h) = \sum_{n\le N}r(n)r(n+h) - 8\left|2^{k+1}-3\right|
\s\left(\frac{h}{2^k}\right)\frac{N}{h},
$$
where $2^k$ is the highest power of 2 dividing $h$. However it is not hard
to see that
$$
g(h) \;:=\;
{(-1)^h8\over h}\sum_{d|h}(-1)^dd = {8\over h}\left|2^{k+1}-3\right|
\s\left(\frac{h}{2^k}\right).\leqno(3.3)
$$
Namely if $k=0$ then $h$ is odd and both expressions in (3.3) reduce to
$8\s(h)/h$. If $k\ge1$, then setting $H = h/2^k$ the identity becomes
$$
\sum_{d|2^kH,(2,H)=1}(-1)^dd = (2^{k+1}-3)\s(H).
$$
But the left-hand side equals
$$\eqalign{
&\sum_{d|H}\left((-1)^dd + (-1)^{2d}2d + \ldots (-1)^{2^kd}2^kd\right)\cr&
= -\s(H) + (2 + 2^2 + \ldots + 2^k)\s(H)\cr&
= (-1 + 2^{k+1} - 2)\s(H) = (2^{k+1} - 3)\s(H).\cr}
$$
\smallskip

We start from (3.6) of [5], writing
$$
\sum_{n\le t}r(n)r(n+h) = g(h)t + E(t,h)\qquad(h^2 \le t \le T^{10}),
$$
where $g(h)$ is given by (3.3).
We recall the definition made in [5], namely
$$
f(t,h) := \left\{-(\sqrt{t+h} - \sqrt{t})^2 +
\frac{3(2t+h)+2\sqrt{t(t+h)}}{16\pi^2\sqrt{t(t+h)}T}\right\}
t^{-3/4}(t+h)^{-3/4}
$$
and note that, for $h^2 \le t \le T^{10}$,
$$
f(t,h) \ll h^2t^{-5/2} + T^{-1}t^{-3/2},\;
{\d f(t,h)\over\d t} \ll h^2t^{-7/2} + T^{-1}t^{-5/2}.
$$
Then, as in [5], we can write
$$
\sum(T) \= \sum\nolimits_1(T) + \sum\nolimits_2(T),
$$
where
$$
\sum\nolimits_1(T) \;:= \;
\sqrt{\pi}T^{5/2}\sum_{h\le T^5}g(h)\int_{h^2}^{T^{10}}
e^{-\pi^2T(\sqrt{t+h}-\sqrt{t})^2}f(t,h)\d t,
$$
$$
\sum\nolimits_2(T) \;:=\; \sqrt{\pi}T^{5/2}\sum_{h\le T^5}\int_{h^2}^{T^{10}}
e^{-\pi^2T(\sqrt{t+h}-\sqrt{t})^2}f(t,h)\d E(t,h).
$$
We can evaluate $\sum_1(T)$ (which provides the main terms in (2.6)
plus an error term which is certainly $\ll \sqrt{T}$), as in [5].
The main task consists of the estimation of $\sum\nolimits_2(T)$,
which contributes to the error term in (2.6). We effect this
by an integration by parts. The integrated terms
will be small, and we are left with the estimation of
$$
T^{5/2}\sum_{h\le T^5}\int_{h^2}^{T^{10}}
E(t,h)u(t,h)\d t,
$$
where ($h^2 \le t \le T^{10}$)
$$
\eqalign{
u(t,h) &= {\d\over\d t}\left(
e^{-\pi^2T(\sqrt{t+h}-\sqrt{t})^2}f(t,h)\right)\cr&
\ll e^{-{2Th^2\over t}}\left(h^2t^{-7/2} + T^{-1}t^{-5/2}
+ Th^4t^{-9/2}\right).\cr}\leqno(3.4)
$$
Now write the above sum as
$$
T^{5/2}\int_{1}^{T^{10}}\sum_{h\le t^{1/2}}
E(t,h)u(t,h)\d t,
$$
and divide the intervals of integration and summation into $O(\log^2T)$
subintervals of the form $[K, 2K]$ and $[H, 2H]$, respectively.

Note that (3.1) can be used with
$$
u(t,m) \ll \a_m =
e^{-{Tm^2\over K}}\left(m^2K^{-7/2} + T^{-1}K^{-5/2}
+ Tm^4K^{-9/2}\right),
$$
since the dependence of $u(t,m)$ on $t$ when $t \in [K,2K]$ is harmless.
Thus we obtain a contribution which is
$$
\eqalign{
& \ll T^{5/2}\log^2T\max_{K\ll T^{10},H\ll\sqrt{K}}
\int_K^{2K}\Big|\sum_{H<h\le2H}E(t,h)u(t,h)\Big|\d t\cr&
\ll_\e T^{5/2+\e}\max_{K\ll T^{10},H\ll\sqrt{K}}
e^{-TH^2/K}KH^{1/2}(K^{2/3}H^{1/2}+K^{1/3}H^{5/6})\times\cr&
\times(H^2K^{-7/2} + T^{-1}K^{-5/2} + TH^4K^{-9/2})\cr&
\ll_\e T^{5/2+\e}\max_{K\ll T^{10},H\ll\sqrt{K}}
e^{-{TH^2\over K}}(H^3K^{-11/6} + HT^{-1}K^{-5/6} + TH^5K^{-17/6}),
\cr}\leqno(3.5)
$$
since $K^{1/3}H^{5/6} \le K^{2/3}H^{1/2}$. Now using
$$
e^{-x}x^\a \;\le\; e^{-\a}\a^\a \;\ll\;1
\qquad(x \ge 0,\;\a > 0 \;{\roman {fixed}})
$$
we obtain
$$
e^{-{TH^2\over K}}H^3K^{-{11\over6}}
 \;\le\; e^{-{TH^2\over K}}\left({TH^2\over K}\right)^{11\over6}
T^{-{11\over6}} \ll T^{-{11\over6}},
$$
and likewise
$$
e^{-{TH^2\over K}}HT^{-1}K^{-{5\over6}} \ll T^{-{11\over6}},\;
e^{-{TH^2\over K}}TH^5K^{-{17\over6}}
\ll T^{-{11\over6}}.
$$
Since ${5\over2} - {11\over6} = {2\over3}$, we obtain
$$
\sum\nolimits_2(T) \;\ll_\e\; T^{{2\over3}+\e},
$$
which gives then (2.6).

Alternatively we may use (3.2) (although this is not uniform in $m$,
it is crucial that the exponent of $m$ is small),
namely $E(N,m) \ll_\e N^{{2\over3}+\e}m^\b$
with $\b = 5/42$. Since $\exp(-TH^2/K) \le T^{-50}$ for $H \ge \left(50K
\log T/T\right)^{1/2}$ and $H\ge1$ has to hold, it follows that in the
first bound in (3.5) it suffices to take the maximum
 over $1 \le H \le \left(50K
\log T/T\right)^{1/2}$ and $T/\log T \ll K \ll T^{50}$. Trivial estimation,
(3.2) and (3.4) yield then a contribution which is
$$\eqalign{
\ll T^{{5\over2}+\e}&\max_{{T\over\log T}\ll K\ll T^{10}}
\left\{{({K\over T})}^{{3+\b}\over2}K^{{2\over3}-{5\over2}}
+ {({K\over T})}^{{1+\b}\over2}T^{-1}K^{{2\over3}-{3\over2}}
+ T{({K\over T})}^{{5+\b}\over2}K^{{2\over3}-{7\over2}}\right\}
\cr&
\ll_\e T^{{2\over3}+\e}\cr}
$$
provided that $0 \le \b \le 2/3$, which in our case is amply satisfied
since we can take $\b = 5/42$. This furnishes another proof
of the Theorem.

\Refs

\bigskip
\item{[1]} F. Chamizo, {\it Correlated sums of $r(n)$}, J. Math. Soc.
Japan  {\bf 51}(1999), 237-252.

\item{[2]} G.H. Hardy, {\it The average order of the arithmetical
functions $P(x)$ and $\D(x)$}, Proc. London Math. Soc. (2){\bf15}
(1916), 192-213.

\item{[3]} M.N. Huxley, {\it Exponential sums and lattice points
II}, Proc. London Math. Soc. (3){\bf66}(1993), 279-301.

\item{[4]} A. Ivi\'c, {\it The Riemann zeta-function}, John Wiley
\& Sons, New York, 1985.

\item{[5]} A. Ivi\'c, {\it The Laplace transform of the square in
the circle and divisor problems}, Studia Scien. Math. Hungarica
{\bf32}(1996), 181-205.

\item{[6]} D. Ismoilov, {\it Additive divisor problems} (in Russian),
Tad\v zik State University, Dushanbe, 1988.

\item{[7]} M. Jutila, Lectures on a method in the theory of
exponential sums, Vol. {\bf80}, Tata Institute of fundamental
research, Bombay, distr. by Springer Verlag, Berlin etc., 1987.

\item{[8]} M. Jutila, {\it Exponential sums connected with
quadratic forms}, in ``Number Theory" (ed. R.A.  Mollin), Walter de
Gruyter,  Berlin etc. 1990,  271-286.

\item{[9]} I. K\'atai, {\it The number of lattice points in
a circle} (in Russian), Ann. Univ. Sci. Budapest E\"otv\"os Sect.
Math. {\bf8}(1965), 39-60.

\item{[10]} T. Meurman, {\it On the binary divisor problem},  in
 ``Number Theory", Proc. Conference in honor of K. Inkeri
(Turku, 1999, eds. M. Jutila and
T. Mets\"ankyl\"a), Walter de Gruyter, Berlin etc., 2001, pp. 223-246.

\item{[11]} Y. Motohashi, {\it The binary additive divisor problem},
Ann. Sci. \'Ecole Normale Sup\'erieure (4){\bf 27}(1994), 529-572.

\bigskip

Aleksandar Ivi\'c

Katedra Matematike RGF-a, Universitet u Beogradu

\DJ u\v sina 7, 11000 Beograd

Serbia and Montenegro, \tt aivic\@rgf.bg.ac.yu


\endRefs
\vfill

\bye